\documentclass{article}    
\usepackage{amstext}
\def\qed{$\rlap{$\sqcap$}\sqcup$}
\usepackage{latexsym}

\begin{document}           

\begin{center}
{\ }\\
{\huge {\bf The Hilbert functions which force\\ the Weak Lefschetz Property}} \\ [.250in]
{JUAN MIGLIORE\\
Department of Mathematics, University of Notre Dame, Notre Dame, IN 46556, USA
\\E-mail: Juan.C.Migliore.1@nd.edu}
{{\ }
\\FABRIZIO ZANELLO\\
Department of Mathematics, Royal Institute of Technology, 100 44 Stockholm, Sweden
\\E-mail: zanello@math.kth.se\\
(Current address: Department of Mathematics, University of Notre Dame, Notre Dame, IN 46556, USA)}
\end{center}

{\ }\\
\\
ABSTRACT. The purpose of this note is to characterize the finite Hilbert functions which force all of their artinian algebras to enjoy the Weak Lefschetz Property (WLP). Curiously, they turn out to be exactly those (characterized by Wiebe in $[Wi]$) whose Gotzmann ideals have the WLP.\\
This implies that, if a Gotzmann ideal has the WLP, then all algebras with the same Hilbert function (and hence lower Betti numbers) have the WLP as well. However, we will answer in the negative, even in the case of level algebras, the most natural question that one might ask after reading the previous sentence: If $A$ is an artinian algebra enjoying the WLP, do all artinian algebras with the same Hilbert function as $A$ and Betti numbers lower than those of $A$ have the WLP as well?\\
Also, as a consequence of our result, we have another (simpler) proof of the fact that all codimension 2 algebras enjoy the WLP (this fact was first proven in $[HMNW]$, where it was shown that even the Strong Lefschetz Property holds).\\

{\large

{\ }\\
\\\indent
Let $A=R/I$ be a standard graded artinian algebras, where $R$ is a polynomial ring in $r$ variables over a field $k$ of characteristic zero, $I$ is a homogeneous ideal of $R$, and the $x_i$'s all have degree 1. We say that $A$ enjoys the {\it Weak Lefschetz Property} ({\it WLP}) if, for a generic linear form $L\in R$ and for all indices $i\geq 0$, the multiplication map \lq \lq $\cdot L$" between the $k$-vector spaces $A_i$ and $A_{i+1}$ has maximal rank (notice that, since $A$ is artinian, $A_i=0$ for $i\gg 0$, and therefore only a finite number of maps have to be considered).\\\indent
The WLP is a fundamental and very natural property of artinian algebras, and has recently received much attention. For a broad overview and the main results achieved so far, see $[GHMS]$, $[HMNW]$ and $[MM]$ (along with their bibliographies), as well as three more recent works of M. Boij and the two authors of this note: $[Za]$, $[Mi]$, $[BZ]$.\\\indent
In particular, one interesting problem that immediately arose was to study the structure of the Hilbert functions of the algebras having the WLP. In $[HMNW]$, the four authors characterized the Hilbert functions forcing {\it at least one} of their artinian algebras to have the Weak Lefschetz Property. See also Harima's paper $[Ha]$, where the Hilbert functions of Gorenstein algebras having the WLP are classified.\\
\\\indent
In this brief note, we will characterize the Hilbert functions which force {\it all} of their artinian algebras to enjoy the WLP. A curious fact is that these Hilbert functions will be exactly those (characterized by Wiebe in $[Wi]$, Theorem 4.2) whose Gotzmann ideals have the WLP.\\\indent
As a consequence of our result, we have another (simpler) proof of the fact that all codimension 2 artinian algebras enjoy the WLP (this was first proven in $[HMNW]$, Proposition 4.4, where it was shown that even the Strong Lefschetz Property holds).\\\indent
Our characterization, along with the result of Wiebe, immediately suggests the following question: If $A$ is an artinian algebra enjoying the WLP, do all artinian algebras with the same Hilbert function as $A$ and Betti numbers lower than those of $A$ have the WLP as well?\\\indent
We will prove that this is not the case, by exhibiting examples of codimension 3 level algebras (which are artinian reductions of algebras of reduced sets of points of the projective space ${\bf P}^3$), say $A_1$ with the WLP and $A_2$ without the WLP, such that the Betti numbers of $A_2$ are obtained from those of $A_1$ by performing cancellations. We will also show that, even in an instance when all level algebras having a given Hilbert function $H$ do not enjoy the WLP (their existence is proven by Boij and the second author in $[BZ]$, Theorem 3.2), we can obtain (of course non-level) algebras having the same Hilbert function $H$ and enjoying the WLP by adding suitable redundant terms to a level minimal free resolution (MFR).\\
\\\indent
Let us now state the main definitions and results we will need in this paper.\\
\\\indent 
{\bf Definition-Remark 1.} Let $n$ and $i$ be positive integers. The {\it i-binomial expansion of n} is $$n_{(i)}={n_i\choose i}+{n_{i-1}\choose i-1}+\cdots +{n_j\choose j},$$ where $n_i>n_{i-1}>...>n_j\geq j\geq 1$. Under these hypotheses, the $i$-binomial expansion of $n$ is unique (e.g., see $[BH]$, Lemma 4.2.6).\\\indent 
Furthermore, define $$(n_{(i)})_{-1}^{-1}={n_i-1\choose i-1}+{n_{i-1}-1\choose i-1-1}+\cdots +{n_j-1\choose j-1},$$
and $$n_{<i>}={n_i-1\choose i}+{n_{i-1}-1\choose i-1}+\cdots +{n_j-1\choose j},$$
where we set ${m\choose q}=0$ whenever $m<q$.\\
\\\indent 
{\bf Theorem 2} (Green). {\it Let $h_d$ be the entry of degree $d$ of the Hilbert function of $R/I$ and let $L$ be a generic linear form of $R$. Then the degree $d$ entry $h_d^{'}$ of the Hilbert function of $R/(I,L)$ satisfies the inequality
$$h_d^{'}\leq (h_d)_{<d>}.$$}
\\\indent
{\bf Proof.} See $[Gr]$, Theorem 1.{\ }{\ }\qed \\
\\\indent
The following is a simple - but fundamental - observation of Stanley (see $[St]$, bottom of p. 67, where this remark appears in a more general form):\\
\\\indent
{\bf Lemma 3} (Stanley). {\it Let $R/I$ be an artinian algebra, and let $L\in R$ be a linear form not belonging to $I$. Then the Hilbert function $H$ of $R/I$ can be written as
$$H: {\ }h_0=c_0=1,h_1=b_1+c_1,...,h_e=b_e+c_e,h_{e+1}=0,$$ where
$$B: {\ }b_1=1,b_2,...,b_e,b_{e+1}=0$$
is the Hilbert function of $R/(I:L)$ (with the indices shifted by 1 to the left), and
$$C: {\ }c_0=1,c_1,...,c_e,c_{e+1}=0$$
is the Hilbert function of $R/(I,L)$.}\\
\\\indent
{\bf Proof.} This decomposition of $H$ is an immediate consequence of the fact that
$$0\longrightarrow R/(I:L) (-1)\stackrel{\cdot L}{\longrightarrow} R/I\stackrel{\cdot 1}{\longrightarrow} R/(I,L)\longrightarrow 0$$
is an exact complex.{\ }{\ }\qed \\
\\\indent
Finally, we recall a theorem of Wiebe, which will be fundamental for us here:\\
\\\indent
{\bf Theorem 4} (Wiebe). {\it Let $A$ be a Gotzmann algebra with Hilbert function\\$H$: $1,h_1,h_2...,h_e, h_{e+1}=0$, and let $t$ be the smallest integer such that $h_t\leq t$. Then $A$ enjoys the WLP if and only if, for all $i=1,2,...,t-1$, we have $$h_{i-1}=((h_i)_{(i)})^{-1}_{-1}.$$}
\\\indent
{\bf Proof.} See $[Wi]$, Theorem 4.2.{\ }{\ }\qed \\
\\\indent
The formula for $H$ provided by Wiebe's characterization basically means that $h_{t-1}$ can be arbitrary (of course subject to Macaulay's theorem) but that, by Bigatti-Geramita's $[BG]$, Lemma 3.3, any other $h_i$, for $i\leq t-2$, is uniquely determined as the minimal possible entry of degree $i$ of $H$, given $h_{i+1}$. Equivalently, the lex-segment (and hence every) algebra having Hilbert function $H$ cannot have socle in degrees lower than $t-1$. However, conversely, at each step the growth from degree $i$ to degree $i+1$ is not necessarily maximal (e.g., $h_2=6$ grows maximally to 10 in degree 3, but $((h_3)_{(3)})^{-1}_{-1}=6$ for all $h_3=8,9,10$).\\
\\\indent
We are now ready to prove our main result, namely the characterization of the finite Hilbert functions which force {\it all} of their (artinian) algebras to enjoy the WLP. Strangely enough, they turn out to be same which force only their {\it Gotzmann} (artinian) algebras to enjoy the WLP!\\
\\\indent
{\bf Theorem 5.} {\it Let $H$: $1,h_1,h_2...,h_e, h_{e+1}=0$ be a possible Hilbert function (according to Macaulay's theorem), and let $t$ be the smallest integer such that $h_t\leq t$. Then all the artinian algebras having Hilbert function $H$ enjoy the WLP if and only if, for all $i=1,2,...,t-1$, we have $$h_{i-1}=((h_i)_{(i)})^{-1}_{-1}.$$}
\\\indent
{\bf Proof.} The \lq \lq only if" part immediately follows from Wiebe's Theorem 4 (if all algebras with Hilbert function $H$ satisfy the WLP, then in particular the Gotzmann algebras do, and hence $H$ has the desired form).\\\indent
Hence assume that $h_{i-1}=((h_i)_{(i)})^{-1}_{-1}$ for $0<i<t$. We want to prove that all artinian algebras $A=R/I$ with Hilbert function $H$ enjoy the WLP. Since, clearly, $H$ increases up to degree $t-1$ and $h_{t-1}\geq h_t$ (by definition of $t$ and by construction of $H$), then by the definition of Weak Lefschetz Property it is easy to see that $A$ has the WLP if and only if, for any generic linear form $L\in R$, the Hilbert function of $R/(I,L)$ is \begin{equation}\label{ee}1,h^{'}_1=h_1-1,h^{'}_2=h_2-h_1,...,h^{'}_{t-1}=h_{t-1}-h_{t-2},h^{'}_t=0.\end{equation}\indent
By Green's Theorem 2, since $h_t\leq t$, we clearly have
$$h^{'}_t\leq (h_t)_{<t>}=\left({t\choose t}+{t-1\choose t-1}+...{t-h_t+1\choose t-h_t+1}\right)_{<t>}=0,$$ whence $h^{'}_t=0$, as we wanted to show.\\\indent
Thus, from now on, let $0<i<t$. It is a trivial exercise, by the Pascal triangle inequality, to show that $h_{i}-(h_i)_{<i>}=((h_i)_{(i)})^{-1}_{-1}$. Hence, by Green's theorem, we have that
\begin{equation}\label{ii}h^{'}_i\leq (h_i)_{<i>}=h_i-((h_i)_{(i)})^{-1}_{-1}=h_i-h_{i-1}.\end{equation}\indent
On the other hand, since $I\subset I:L$, by Lemma 3 we must have that $h_i-h^{'}_i\leq h_{i-1}$, i.e. that
\begin{equation}\label{iii}h^{'}_i\geq h_i-h_{i-1}.\end{equation}\indent
Thus, from (\ref{ii}) and (\ref{iii}), we get (\ref{ee}), and the proof of the theorem is complete.{\ }{\ }\qed \\
\\\indent
As we noticed before, by Wiebe's theorem, our Theorem 5 can be rephrased in the following surprising form:\\
\\\indent
{\bf Corollary 6.} {\it Let $H$ be a finite Hilbert function. Then all algebras having Hilbert function $H$ enjoy the WLP if and only if {\em any Gotzmann} algebra having Hilbert function $H$ enjoys the WLP.}\\
\\\indent
A very interesting consequence to our Theorem 5 is the following result (that, as we said in the introduction, was first proven in $[HMNW]$, Proposition 4.4):\\
\\\indent
{\bf Corollary 7} ($[HMNW]$). {\it All codimension 2 artinian algebras enjoy the WLP.}\\
\\\indent
{\bf Proof.} It is immediate to see, by the growth condition imposed by Macaulay's theorem on Hilbert functions $H$ starting with 1,2,..., that they must all have the form $H: {\ } h_0=1,h_1=2,...,h_{t-1}=t,h_t\leq t,...$.\\\indent
Since, clearly, $((i+1)_{(i)})^{-1}_{-1}=i$, the result immediately follows from Theorem 5.{\ }{\ }\qed \\
\\\indent
Also notice that, from our Theorem 5, one implication of Wiebe's Theorem 4 immediately follows (of course, not the same that we used in our proof!):\\
\\\indent
{\bf Corollary 8} (Wiebe). {\it Let $H$: $1,h_1,h_2...,h_e, h_{e+1}=0$ be the Hilbert function of an artinian algebra such that, for all $i=1,2,...,t-1$, $h_{i-1}=((h_i)_{(i)})^{-1}_{-1}.$ Then all Gotzmann ideals with Hilbert function $H$ satisfy the WLP.}\\
\\\indent 
Since the Betti numbers of lex-segment ideals - and therefore, more generally, of Gotzmann ideals - are the largest possible for all the algebras with a given Hilbert function (see $[Bi]$, $[Hu]$ and $[Pa]$), at this point it is natural to wonder if our theorem is the consequence of a more general fact: that is, is it true that, given an artinian algebra $A_1$ with the Weak Lefschetz Property, then all the algebras $A_2$, with the same Hilbert function as $A_1$ but Betti numbers strictly lower than those of $A_1$, enjoy the WLP as well?\\\indent
Unfortunately (or not?), the answer is negative, as we will show in the next proposition. Indeed, we can even provide counterexamples which are constructed as artinian reductions of one-dimensional algebras of points of ${\bf P}^3$, where both $A_1$ and $A_2$ are level algebras.\\\indent
Also, we will show that even in cases (whose existence was proven in $[BZ]$, Theorem 3.2) where a Hilbert function $H$ forces all of its level algebras to not enjoy the WLP, one can find (of course non-level) algebras, with the same Hilbert function $H$ and Betti numbers obtained by adding redundant terms to a level MFR, which enjoy the WLP.\\\indent
Notice that, when $A_1$ and $A_2$ have exactly the same Betti numbers, the fact that the answer to the above question is negative is already known, even in the Gorenstein case - see $[HMNW]$, Example 3.10, where a codimension 4 Gorenstein algebra $A_2$ with the WLP is provided, with $A_2$ having the same Betti numbers and the same Hilbert function 1,4,10,10,4,1,0 of the well-known non-WLP example of Ikeda ($[Ik]$, Example 4.4).\\
\\\indent
{\bf Proposition 9.} {\it i). There exist reduced sets of points $Y_1, Y_2$ in ${\bf P}^3$ with the same Hilbert function, and both having {\em level} coordinate ring, such that their (arbitrary) artinian reductions, $A_1$ and $A_2$ respectively, are such that $A_1$ has the WLP while $A_2$ does not, and the graded Betti numbers of $A_2$ can be obtained from those of $A_1$ by cancellation of redundant terms.\\\indent
ii). There exists a finite Hilbert function which forces all of its level algebras to fail to have the WLP, but such that one of its non-level algebras has Betti numbers larger than those of a level algebra and enjoys the WLP.}\\
\\\indent
{\bf Proof}. i). Our artinian level reductions $A_1$ and $A_2$, quotients of $R=k[x,y,z]$, will both have Hilbert function
\begin{equation} \label{hf}
1,3,6,10,12,12,0.\end{equation}\indent
We begin with $A_1$.  In this case it is somewhat more intuitive to first construct the artinian algebra and then lift to points $Y_1$.  The idea is based on the fact that (\ref{hf}) coincides with the Hilbert function of 12 points in ${\bf P}^2$ through degree 5, and such an ideal has projective dimension 1.\\\indent
Let $Z$ be a $k$-configuration of type (1,2,4,5); $Z$ consists of 12 points in ${\bf P}^2$ with Hilbert function 
\begin{equation} \label{pointsA}
1,3,6,10,12,12,12,\dots .\end{equation}\indent
Specifically, $Z$ consists of a union of five points on a line $l_1$, four points on a line $l_2$, two points on a line $l_3$ and one additional point.  By $[GHS]$, $Z$ has extremal Betti numbers among sets of points in ${\bf P}^2$ with Hilbert function (\ref{hf}). The Betti diagram for $R/I_Z$ is

\begin{center}

\begin{tabular}{ccccccccccccccccc}
1 & - & - \\
- & - & -  \\
- & - & - \\
- & 3 & 2 \\
- & 2 & 2
\end{tabular}

\end{center}

\indent Now we consider the ideal $I = I_Z + (x,y,z)^6$.  Clearly $R/I$ has Hilbert function (\ref{hf}), and the following level Betti diagram is forced from the Hilbert function and the above Betti diagram:

\begin{center}

\begin{equation} \label{betti A1}
\begin{tabular}{ccccccccccccccccc}
1 & - & - & - \\
- & - & -  & - \\
- & - & - & - \\
- & 3 & 2 & - \\
- & 2 & 2 & - \\
- & 12 & 24 & 12
\end{tabular}
\end{equation}

\end{center}

\indent Finally, $A_1 = R/I$ clearly has the WLP since up to and including degree 5 it coincides with the coordinate ring of a set of points in ${\bf P}^2$, which has depth 1.  To construct the set of points, $Y_1$, described in the statement of the proposition, let $C$ be a corresponding arithmetically Cohen-Macaulay union of plane curves of total degree 12 (for instance take $C$ to be a cone over $Z$).  By $[GMR]$, there exists a set, $Y_1$, of 44 points on $C$ whose Hilbert function is the desired truncation.\\\indent
We now construct the level algebra $A_2$ using the methods of $[Mi]$. We will find a set of points $Y_2$ in ${\bf P}^3$ directly, as the union $Y_2 = X_1 \cup X_2 $, based on the computation (in the notation of $[Mi]$):

\begin{equation} \label{calculate}
\begin{tabular}{l|ccccccccccccc}
deg & 0 & 1 & 2 & 3 & 4 & 5 \\ \hline
$\Delta h_{X_1}$ & & 1 & 3 & 6 & 7 & 6 \\
$\Delta h_{X_2}$ & 1 & 2 & 3 & 4 & 5 & 6 \\ \hline
$\Delta h_{Y_2}$ & 1 & 3 & 6 & 10 & 12 & 12
\end{tabular}
\end{equation}

\indent To achieve this, let $X$ be a general set of 4 points in ${\bf P}^3$ and let $W$ be a complete intersection of type (3,3,3) containing $X$. Letting $X_1$ be the residual to $X$ in $W$, clearly $X_1$ has $h$-vector given by the first line of (\ref{calculate}). The MFR of $I_{X_1}$ is easily seen to be 
\begin{equation} \label{res x1}
0 \rightarrow R(-7)^6 \rightarrow R(-6)^{11} \rightarrow 
\begin{array}{c}
R(-3)^3 \\
\oplus \\
R(-5)^3
\end{array}
\rightarrow I_{X_1} \rightarrow 0.
\end{equation}
\indent It is not hard to see that one can construct an arithmetically Cohen-Macaulay curve $C$ of degree 21, with $h$-vector given by the second line of (\ref{calculate}), that contains a complete intersection of type (3,3,3).  Indeed, one can construct $C$ having as a component a complete intersection $C_1$ of type (3,3), and cut $C_1$ by a cubic hypersurface.  Let $X_2$ be a general hyperplane section of $C$.  Then $X_2$ also has $h$-vector given by the second line of (\ref{calculate}), and we take $Y_2 = X_1 \cup X_2$.  This is a basic double G-link, in the language of $[KMMNP]$, and the homogeneous ideal of $Y_2$ has the form
\[
I_{Y_2} = L \cdot I_{X_1} + I_C,
\]
where $L$ is the linear form defining the hyperplane that cuts $X_2$. In particular, $Y_2$ has $h$-vector (\ref{hf}), as desired.\\\indent
Since $I_C$ clearly has MFR
\[
0 \rightarrow R(-7)^6 \rightarrow R(-6)^7 \rightarrow I_C \rightarrow 0,
\]
we use a mapping cone over the short exact sequence
\begin{equation} \label{res x2}
0 \rightarrow I_C(-1) \rightarrow I_{X_1}(-1) \oplus I_C \rightarrow I_{Y_2} \rightarrow 0
\end{equation}
and the MFR's (\ref{res x1}) and (\ref{res x2}) to obtain the following level Betti diagram for $R/I_{Y_2}$:

\begin{center}

\begin{equation} \label{betti A2}
\begin{tabular}{ccccccccccccccccc}
1 & - & - & - \\
- & - & -  & - \\
- & - & - & - \\
- & 3 & - & - \\
- & - & - & - \\
- & 10 & 24 & 12
\end{tabular}
\end{equation}

\end{center}

\indent Notice that there is no possible cancellation to worry about in the mapping cone.  Also, there are clearly no redundant terms.  It is not hard to show, using $[Mi]$, Proposition 2.3, that any artinian reduction of $Y_2$ fails to have the WLP. (The key point is that for a general line $\lambda$ in ${\bf P}^3$, the point of intersection of $\lambda$ with $H_4$ fails to impose a condition on the linear system of quintics containing $Y_2$.)  Finally, note that the Betti numbers of (\ref{betti A2}) can be obtained from those of (\ref{betti A1}) by cancellation (and that this involves only the first two free modules of the MFR, since the level property is preserved). This proves point i).\\\indent
ii). We will consider the case $e=9$ of $[BZ]$, Theorem 3.2, where the Hilbert functions provided are such that all of their level algebras do not enjoy the WLP.  This means that we are considering the Hilbert function
\begin{equation} \label{bzhf}
1,3,5,7,9,11,11,8,5,2,0.
\end{equation}\indent
As above, let $R = k[x,y,z]$. The Hilbert function (\ref{bzhf}) satisfies the condition of $[HMNW]$ for algebras having the WLP, and the construction there gives an algebra with the WLP, and having  maximal Betti numbers {\it among algebras with the WLP} having the given Hilbert function.  In this case one can check (with $[CoCoA]$) that the algebra $R/I_1$ constructed has MFR 
\begin{equation} \label{max wlp betti}
0 \rightarrow 
\begin{array}{c}
R(-9)^3 \\
\oplus \\
R(-10)^3 \\
\oplus \\
R(-11)^3 \\
\oplus \\
R(-12)^2 
\end{array}
\rightarrow 
\begin{array}{c}
R(-7)^2 \\
\oplus \\
R(-8)^6 \\
\oplus \\
R(-9)^6 \\
\oplus \\
R(-10)^6 \\
\oplus \\
R(-11)^4
\end{array}
\rightarrow 
\begin{array}{c}
R(-2) \\
\oplus \\
R(-6)^2 \\
\oplus \\
R(-7)^3 \\
\oplus \\
R(-8)^3 \\
\oplus \\
R(-9)^3 \\
\oplus \\
R(-10)^2
\end{array}
\rightarrow R \rightarrow R/I_1 \rightarrow 0.
\end{equation}

\indent On the other hand, the proof of $[BZ]$, Theorem 3.2 provides a level algebra, $R/I_2$, that can be checked to have MFR
\[
0 \rightarrow R(-12)^2 \rightarrow 
\begin{array}{c}
R(-7)^2 \\
\oplus \\
R(-8)^3 \\
\oplus \\
R(-10) \\
\oplus \\
R(-11)
\end{array}
\rightarrow 
\begin{array}{c}
R(-2) \\
\oplus \\
R(-6)^2 \\
\oplus \\
R(-7)^3 
\end{array}
\rightarrow R \rightarrow R/I_2 \rightarrow 0
\]
\indent This algebra clearly has Betti numbers (much) below the Betti numbers of $R/I_1$, as we wanted to show.{\ }{\ }\qed \\
\\\indent
{\bf Remark 10.} i). By Wiebe's theorem, it is immediate to see that a Gotzmann algebra with Hilbert function (\ref{bzhf}) does not enjoy the WLP. Therefore, in the proof of Proposition 9, we have actually shown that there exist algebras failing to have the WLP whose Betti numbers are above those of an algebra, $R/I_1$, which has the WLP, as well as other algebras, again without the WLP, whose Betti numbers are below those of $R/I_1$.\\\indent
ii). Notice that, actually, for {\it all} the Hilbert functions $H$ provided by $[BZ]$, Theorem 3.2, which force all of its level algebras to not have the WLP, the same argument given in the proof of point ii) of Proposition 9 (employing the condition for algebras having the WLP given in $[HMNW]$) can be used to show the existence of non-level algebras with Hilbert function $H$ enjoying the WLP.\\
\\\indent
{\bf Acknowledgements.} The second author is funded by the G\"oran Gustafsson Foundation.\\
\\
\\
\\
\\
\\
{\bf \huge References}\\
\\
$[Bi]$ {\ } A. Bigatti: {\it Upper bounds for the Betti numbers of a given Hilbert function}, Comm. in Algebra 21 (1993), No. 7, 2317-2334.\\
$[BG]$ {\ } A.M. Bigatti and A.V. Geramita: {\it Level Algebras, Lex Segments and Minimal Hilbert Functions}, Comm. in Algebra 31 (2003), 1427-1451.\\
$[BZ]$ {\ } M. Boij and F. Zanello: {\it Level Algebras with Bad Properties}, Proc. of the Amer. Math. Soc., to appear (preprint: math.AC/0512198).\\
$[BH]$ {\ } W. Bruns and J. Herzog: {\it Cohen-Macaulay rings}, Cambridge studies in advanced mathematics, No. 39, Revised edition (1998), Cambridge, U.K..\\
$[CoCoA]$ A. Capani, G. Niesi and L. Robbiano: {\it CoCoA, a system for doing computations in commutative algebra}, available via anonymous ftp, cocoa.dima.unige.it.\\
$[GHMS]$ {\ } A.V. Geramita, T. Harima, J. Migliore and Y.S. Shin: {\it The Hilbert Function of a Level Algebra}, Memoirs of the Amer. Math. Soc., to appear.\\
$[GHS]$ {\ } A.V. Geramita, T. Harima and Y.S. Shin: {\it Extremal point sets and Gorenstein ideals}, Adv. Math. 152 (2000),  No. 1, 78-119.\\
$[GMR]$ {\ } A.V. Geramita, P. Maroscia and L. Roberts: {\it The Hilbert Function of a Reduced $k$-Algebra}, J. London Math. Soc. 28 (1983), 443-452.\\
$[Gr]$ {\ } M. Green: {\it Restrictions of linear series to hyperplanes, and some results of Macaulay and Gotzmann}, Algebraic curves and projective geometry (1988), 76-86, Trento; Lecture Notes in Math. 1389 (1989), Springer, Berlin.\\
$[Ha]$ {\ } T. Harima: {\it Characterization of Hilbert functions of Gorenstein Artin algebras with the weak Stanley property}, Proc. of the Amer. Math. Soc. 123 (1995), No. 12, 3631-3638.\\
$[HMNW]$ {\ } T. Harima, J. Migliore, U. Nagel and J. Watanabe: {\it The Weak and Strong Lefschetz Properties for Artinian $K$-Algebras}, J. of Algebra 262 (2003), 99-126.\\
$[Hu]$ {\ } H. Hulett: {\it Maximum Betti numbers of homogeneous ideals with a given Hilbert function}, Comm. in Algebra 21 (1993), No. 7, 2335-2350.\\
$[Ik]$ {\ } H. Ikeda: {\it Results on Dilworth and Rees numbers of Artinian local rings}, Japan. J. of Math. 22 (1996), 147-158.\\
$[KMMNP]$ {\ } J. Kleppe, J. Migliore, R.M. Mir\'o-Roig, U. Nagel and C. Peterson: {\it Gorenstein Liaison, Complete Intersection Liaison Invariants and
Unobstructedness}, Memoirs of the Amer. Math. Soc. 154 (2001), 116 pp., Softcover, ISBN 0-8218-2738-3.\\
$[Mi]$ {\ } J. Migliore: {\it The geometry of the Weak Lefschetz Property}, Canadian J. of Math., to appear (preprint: math.AC/0508067).\\
$[MM]$ {\ } J. Migliore and R. Mir\'o-Roig: {\it Ideals of general forms and the ubiquity of the Weak Lefschetz property},
J. of Pure and Applied Algebra 182 (2003), 79-107.\\
$[MN]$ {\ } J. Migliore and U. Nagel: {\it Monomial Ideals and the Gorenstein Liaison Class of a Complete Intersection}, Compositio Math. 133 (2002), 25-36.\\
$[Pa]$ {\ } K. Pardue: {\it Deformation classes of graded modules and maximal Betti numbers}, Illinois J. Math. 40 (1996), 564-585.\\
$[St]$ {\ } R. Stanley: {\it Combinatorics and Commutative Algebra}, Second Ed., Progress in Mathematics 41 (1996), Birkh\"auser, Boston, MA.\\
$[Wi]$ {\ } A. Wiebe: {\it The Lefschetz Property for componentwise linear ideals and Gotzmann ideals}, Comm. in Algebra 32 (2004), No. 12, 4601-4611.\\ 
$[Za]$ {\ } F. Zanello: {\it A non-unimodal codimension 3 level $h$-vector}, J. of Algebra 305 (2006), No. 2, 949-956.\\

}

\end{document}